\documentclass[11pt]{article}
\usepackage{dsfont}
\usepackage{amsmath,amssymb}
\usepackage{amsfonts}
\usepackage{hyperref}
\usepackage{amsthm}
\usepackage{graphicx}
\usepackage{subfigure}
\usepackage{xcolor}
\usepackage{appendix}
\usepackage{cite}
\usepackage{tikz}

\usetikzlibrary{calc,shapes,arrows}

\usepackage{pgfplots}

\usetikzlibrary{patterns}

\usepgfplotslibrary{dateplot}

\renewcommand{\qed}{\hfill\small{$\square$}\normalsize}

\hfuzz=\maxdimen
\theoremstyle{definition}
\newtheorem{lemma}{Lemma}[section]
\newtheorem{definition}[lemma]{Definition}
\newtheorem{proposition}[lemma]{Proposition}
\newtheorem{theorem}[lemma]{Theorem}
\newtheorem{corollary}[lemma]{Corollary}

\newtheorem{remark}{Remark}

\numberwithin{equation}{section}
\renewcommand{\proof}{\textbf{Proof. }}
\renewcommand{\qed}{\hfill\small{$\square$}\normalsize}

\DeclareFixedFont{\Acknowledgment}{OT1}{cmr}{bx}{n}{14pt}
\begin{document}
\title{\bf Parameterized discrete uniformization theorems and curvature flows for polyhedral surfaces, I}
\author{Xu Xu}
\maketitle


\begin{abstract}
In this paper, we introduce a parameterized discrete curvature ($\alpha$-curvature) on polyhedral surfaces,
which is a generalization of the classical discrete curvature.
A discrete uniformization theorem is established for the parameterized discrete curvature, which generalizes
the discrete uniformization theorem obtained by Gu-Luo-Sun-Wu \cite{GLSW}.
We also prove the global rigidity of parameterized discrete curvature with respect to the discrete conformal factors,
which confirms a generalized Luo conjecture \cite{L1} on rigidity of discrete curvatures.
We further introduce a parameterized discrete Yamabe flow for piecewise linear metrics on surfaces.
To handle the possible singularities along the flow, we do surgery on the flow by flipping.
Then we prove that the flow with surgery converges to a piecewise linear metric with constant discrete $\alpha$-curvature,
which confirms another generalized Luo conjecture \cite{L1} on convergence of discrete Yamabe flow with surgery.
We also introduce a parameterized discrete Calabi flow and prove the convergence of the flow with surgery, which generalizes
the convergence result proved in \cite{ZX}.
\end{abstract}



\section{Introduction}\label{section 1}
To study the conformal geometry of piecewise linear metrics on manifolds,
Luo \cite{L1} and R\v{o}cek-Williams \cite{RW} independently introduced
a discrete conformality for piecewise linear metrics (Euclidean polyhedral metrics),
which is now called vertex scaling.
Luo \cite{L1} further introduced the combinatorial Yamabe flow for piecewise linear metrics on triangulated surfaces
and obtained the combinatorial obstruction for the existence of  piecewise linear metrics with constant combinatorial curvature.
Bobenko-Pinkall-Springborn \cite{BPS} proved the global rigidity of vertex scaling
and obtained the relationship between the vertex scaling and the geometry of ideal polyhedra in hyperbolic three space.
They further introduced vertex scaling for piecewise hyperbolic metrics on triangulated surfaces.
Based on Bobenko-Pinkall-Springborn's work \cite{BPS} and Penner's work \cite{Penner},
Gu-Luo-Sun-Wu \cite{GLSW} recently proved a discrete uniformization theorem for piecewise linear metrics
on surfaces via a variational principle established by Luo in \cite{L1}.
Similar discrete uniformization theorem was established by Gu-Guo-Luo-Sun-Wu \cite{GGLSW} for piecewise hyperbolic metrics on surfaces.
Combinatorial Yamabe flow with surgery for polyhedral metrics were defined in \cite{GLSW, GGLSW},
where the long time existence and convergence of the combinatorial
Yamabe flow with surgery were proved.
Following Luo's approach, Ge \cite{Ge-thesis, Ge} introduced the combinatorial Calabi flow on surfaces.
Recently, Zhu and the author \cite{ZX} proved the long-time existence and
convergence of the combinatorial Calabi flow with surgery for vertex scaling of piecewise linear and piecewise hyperbolic metrics on surfaces.
Other related work on vertex scaling could be found in \cite{ GLW, Guo2, LW, LX, Sp, SWGL, Wu, WGS, WX, WZ, XZ AJM}.

In this paper, we introduce a parameterized discrete curvature (combinatorial $\alpha$-curvature) for
piecewise linear metrics with respect to the vertex scaling on surfaces, which is a generalization of the classical discrete curvature.
We prove the global rigidity and a discrete uniformization theorem for this curvature.
We also study the properties of the corresponding discrete curvature flows.
The combinatorial $\alpha$-curvature of piecewise hyperbolic metrics on surfaces is studied in \cite{X3}.
Combinatorial $\alpha$-curvature was introduced by Ge and the author \cite{GX3,GX4}
for Thurston's circle packing metrics as a generalization of the classical
combinatorial curvature.
There are lots of works on combinatorial curvatures and combinatorial curvature flows on surfaces and 3-dimensional manifolds, see \cite{CL,CR,DV,FGH, FGHX, Ge-thesis,Ge,GH,GH1,GJ1,GJ2,GJ3,GJS,GM,GX1,GX2,GXcv,GX3,GX4,GX5,G1,G2,G3,G5,GT,Guo,L2,MR,T1,X1,
X2,X2 JFA2021, X2 JFA2022, X3,X4,Z} and others.

Suppose $S$ is a closed connected surface and $V$ is a finite subset of $S$, $(S, V)$ is called a marked surface.
A piecewise linear metric (PL metric) on $(S, V)$ is a flat cone metric with cone points contained in $V$.
Suppose $\mathcal{T}=(V, E, F)$ is a triangulation of $(S, V)$, where $V,E,F$ represent the set of vertices, edges and faces respectively.
We use $(S, V, \mathcal{T})$ to denote a triangulated surface.
If a map $d: E\rightarrow (0, +\infty)$ satisfies that $d_{rs}<d_{rt}+d_{st}$ for $\{r,s,t\}=\{i,j,k\}$,
where $d_{rs}=d(\{rs\})$ and $\{i,j,k\}$ is any triangle in $F$,
then $d$ determines a PL metric on $(S, V)$, which is still denoted by $d$.
Given $(S, V)$ with a triangulation $\mathcal{T}$ and a map $d: E\rightarrow (0, +\infty)$ determined by a PL metric $d$ on $(S, V)$,
the vertex scaling \cite{L1,RW} of the PL metric $d$
by a function $w: V\rightarrow (0, +\infty)$ is defined to be the PL metric $w*d$ on $(S, V)$
determined by the map $w*d: E\rightarrow (0, +\infty)$  with
\begin{equation*}
(w*d)_{ij}:=w_iw_jd_{ij}, \ \ \forall \{ij\}\in E.
\end{equation*}
The function $w: V\rightarrow (0, +\infty)$ is called a \emph{conformal factor}.
For a triangulated surface $(S, V, \mathcal{T})$ with a PL metric $d$,
we denote the admissible
space of conformal factors by $\Omega^{\mathcal{T}}(d)$, which is the set of discrete conformal factors such that
the triangle inequalities are satisfied for every face in $\mathcal{T}$.
Set $u_i=\ln w_i$, $i=1,\cdots, N$, and
$\mathcal{U}^\mathcal{T}(d)=\ln \Omega^{\mathcal{T}}(d)$. Here $N=|V|$.

Suppose $(S, V)$ is a marked surface with a PL metric $d$.
The combinatorial curvature $K_i$ of $d$ at $v_i\in V$ is $2\pi$ less the cone angle at $v_i$.
If $\mathcal{T}$ is a geometric triangulation of $(S, V)$ with a PL metric $d$, we denote
$\theta_i^{jk}$ as the inner angle at the vertex $v_i$ of the triangle $\triangle ijk$.
Then the combinatorial curvature $K_i=2\pi-\sum_{\triangle ijk\in F}\theta_i^{jk}$.
Note that
the combinatorial curvature $K$ is independent of the geometric triangulations of $(S, V)$ with a PL metric $d$.

\begin{definition}
Suppose $(S, V, \mathcal{T})$ is a triangulated surface with a PL metric $d$ and
$w: V\rightarrow (0, +\infty)$ is a conformal factor of $d$ on $(S, V, \mathcal{T})$.
For any $\alpha\in \mathbb{R}$, the combinatorial $\alpha$-curvature of $w*d$ on $(S, V, \mathcal{T})$ is defined to be
\begin{equation}
R_{\alpha, i}=\frac{K_i}{w_i^\alpha}.
\end{equation}
\end{definition}

In the case that $\alpha=0$, the curvature $R_0$ is the classical combinatorial curvature $K$.
Furthermore, for any constant $\lambda>0$, we have
$R_{\alpha, i}(\lambda*l)=\lambda^{-\alpha}R_{\alpha, i}(l).$
Especially, for $\alpha=1$, we have
$R_{1, i}(\lambda*l)=\lambda^{-1}R_{1, i}(l),$
which is parallelling to the transformation of smooth Gaussian curvature $K_{\lambda g}=\lambda^{-1}K_g$.

For combinatorial $\alpha$-curvature of PL metrics on triangulated surfaces, we have the following global rigidity, which confirms a
generalized Luo conjecture \cite{L1} on rigidity of discrete curvatures.
\begin{theorem}\label{main rigidity intro}
Suppose $(S, V, \mathcal{T})$ is a triangulated closed surface with a PL metric $d$
and $\alpha\in \mathbb{R}$ is a constant such that $\alpha\chi(S)\leq0$.
$\overline{R}$ is a given function defined on $V$.
\begin{description}
  \item[(1)]
  If $\alpha \overline{R}\equiv0$, then there exists at most one conformal factor $\overline{w}\in \Omega^\mathcal{T}(d)$
  with  $\alpha$-curvature $\overline{R}$  up to scaling;
  \item[(2)]
  If $\alpha\overline{R}\leq 0$ and $\alpha\overline{R}\not\equiv 0$, then there
  exists at most one conformal factor $\overline{w}\in \Omega^\mathcal{T}(d)$
  with  $\alpha$-curvature $\overline{R}$.
\end{description}
\end{theorem}

If $\alpha=0$, there is no restriction on $R_0=K$ and the global rigidity of $R_0=K$ in Theorem \ref{main rigidity intro}
is reduced to the rigidity proved in \cite{L1, BPS}.

For $\alpha$-curvature, it is interesting to consider the following Yamabe problem.\\
~\\
\textbf{Combinatorial $\alpha$-Yamabe Problem:} \emph{Suppose $(S, V, \mathcal{T})$ is a closed triangulated surface with a PL metric $d$,
does there exist any conformal factor $w: V\rightarrow (0, +\infty)$ in $\Omega^{\mathcal{T}}(d)$
such that $w*d$ has constant $\alpha$-curvature?}
~\\

In this paper, we prove the following result on combinatorial $\alpha$-Yamabe problem,
which is also called a parameterized discrete uniformization theorem.

\begin{theorem}\label{main theorem existence of constant alpha curvature metric}
Suppose $(S, V)$ is a closed connected marked surface with a PL metric $d_0$ and
$\alpha\in \mathbb{R}$ is a constant such that $\alpha\chi(S)\leq 0$.
Then there exists a PL metric in the conformal class $\mathcal{D}(d_0)$ with constant $\alpha$-curvature.
\end{theorem}

Here $\mathcal{D}(d_0)$ is the discrete conformal class defined in the sense of Gu-Luo-Sun-Wu \cite{GLSW}.
Please refer to Definition \ref{GLSW's definition of edcf} for this.
Theorem \ref{main theorem existence of constant alpha curvature metric} is a parameterized generalization of Gu-Luo-Sun-Wu's discrete uniformization theorem in \cite{GLSW}. Especially, if $\alpha=0$, Theorem \ref{main theorem existence of constant alpha curvature metric} is the discrete uniformization theorem
proved in \cite{GLSW}.


To study the combinatorial $\alpha$-Yamabe problem, we introduce the combinatorial $\alpha$-Yamabe flow and combinatorial
$\alpha$-Calabi flow for PL metrics on surfaces.


\begin{definition}
Suppose $(S, V, \mathcal{T})$ is a triangulated closed surface with a PL metric $d_0$
and $\alpha\in \mathbb{R}$ is a constant.
The normalized combinatorial $\alpha$-Yamabe flow is defined to be
\begin{equation}\label{alpha Yamabe flow}
\left\{
  \begin{array}{ll}
    \frac{dw_i}{dt}=(R_{\alpha, av}-R_{\alpha, i})w_i, & \hbox{  } \\
    w_i(0)=1, & \hbox{  }
  \end{array}
\right.
\end{equation}
where
$R_{\alpha, av}=\frac{2\pi\chi(M)}{\sum_{i=1}^Nw_i^\alpha}.$
\end{definition}
When $\alpha=0$, this is the combinatorial Yamabe flow introduced by Luo \cite{L1}.
By direct calculations,
the combinatorial $\alpha$-curvature $R_\alpha$ evolves according to
\begin{equation}\label{evolution of curv along alpha Yamabe flow}
\frac{dR_{\alpha, i}}{dt}=(\Delta^\mathcal{T}_{\alpha}R_{\alpha})_i+\alpha R_{\alpha, i}(R_{\alpha, i}-R_{\alpha, av})
\end{equation}
along the combinatorial $\alpha$-Yamabe flow (\ref{alpha Yamabe flow}),
where the $\alpha$-Laplace operator $\Delta^\mathcal{T}_{\alpha}$ on $(S, V, \mathcal{T})$ is defined to be
\begin{equation*}
(\Delta^\mathcal{T}_{\alpha}f)_i
=\frac{1}{w_i^\alpha}\sum_{j\sim i} \left(\cot \theta_{k}^{ij}+\cot\theta_l^{ij}\right) (f_j-f_i)
\end{equation*}
for $f\in \mathbb{R}^V$.
Here $\theta_k^{ij}$ and $\theta_l^{ij}$ are two inner angles facing the edge $\{ij\}$.
(\ref{evolution of curv along alpha Yamabe flow}) is similar to the evolution of Gaussian curvature
along the normalized Ricci flow on surfaces \cite{CH1,H2}.

\begin{definition}
Suppose $(S, V, \mathcal{T})$ is a closed triangulated surface with a PL metric $d_0$
and $\alpha\in \mathbb{R}$ is a constant.
The combinatorial $\alpha$-Calabi flow is defined to be
\begin{equation}\label{alpha Calabi flow}
\left\{
  \begin{array}{ll}
    \frac{dw_i}{dt}=(\Delta^\mathcal{T}_\alpha R_{\alpha})_iw_i, & \hbox{ } \\
    w_i(0)=1. & \hbox{ }
  \end{array}
\right.
\end{equation}
\end{definition}
When $\alpha=0$, this is the combinatorial Calabi flow introduced by Ge \cite{Ge-thesis}.
By direct calculations, the combinatorial $\alpha$-curvature $R_\alpha$ evolves according to
\begin{equation*}
\frac{dR_{\alpha, i}}{dt}=-(\Delta^\mathcal{T}_{\alpha})^2R_{\alpha, i}-\alpha R_{\alpha, i}(\Delta^\mathcal{T}_{\alpha}R_\alpha)_i
\end{equation*}
along the combinatorial $\alpha$-Calabi flow (\ref{alpha Calabi flow}),
which is similar to the evolution of Gaussian curvature along the surface Calabi flow \cite{Ca, Ca2, CP}.
If the parameters are chosen properly, the evolution equations of combinatorial $\alpha$-curvature along the combinatorial $\alpha$-Yamabe flow
and $\alpha$-Calabi flow are formally the same as the evolution equations of Gaussian curvature along
the surface Ricci flow and surface Calabi flow respectively. Please refer to \cite{GX4} for this.

The combinatorial $\alpha$-flows (combinatorial $\alpha$-Yamabe flow and combinatorial $\alpha$-Calabi flow)
may develop singularities.
To handle the possible singularities along the combinatorial $\alpha$-flows,
we do surgery on the flows by flipping, the idea of which comes form \cite{GLSW, GGLSW, L1}.
Note that the weight in the $\alpha$-Laplace operator is
$\omega_{ij}=\frac{\cot \theta_{k}^{ij}+\cot\theta_l^{ij}}{w_i^\alpha}$, which is not symmetric with the indices $i$ and $j$.
To ensure that the discrete $\alpha$-Laplace operator have good properties along the $\alpha$-flows,
especially the discrete maximal principle could be applied on the combinatorial $\alpha$-Yamabe flow,
we need the weight $\omega_{ij}$ to be nonnegative on every edge,
which is equivalent to $\theta_{k}^{ij}+\theta_l^{ij}\leq \pi$ for every edge $\{ij\}\in E$.
This is exactly the Delaunay condition on the triangulation \cite{BS}.
This condition is imposed on both the combinatorial $\alpha$-Yamabe flow and the combinatorial $\alpha$-Calabi flow.
Note that every PL metric on $(S, V)$ admits at least one Delaunay triangulation \cite{AK, BS, R},
so this additional condition is reasonable.
Along the $\alpha$-flows on $(S, V)$ with a triangulation $\mathcal{T}$, if $\mathcal{T}$ is Delaunay in $w(t)*d_0$  for
$t\in [0, T]$ and not Delaunay in $w(t)*d_0$  for $t\in (T, T+\epsilon)$, $\epsilon>0$, there exists
an edge $\{ij\}\in E$ such that $\theta_{k}^{ij}(t)+\theta_l^{ij}(t)\leq \pi$
for $t\in [0, T]$ and $\theta_{k}^{ij}(t)+\theta_l^{ij}(t)> \pi$ for $t\in (T, T+\epsilon)$.
Then we replace the triangulation $\mathcal{T}$ by a new triangulation $\mathcal{T}'$ at time $t=T$
via replacing two triangles $\triangle ijk$ and $\triangle ijl$
adjacent to $\{ij\}$ by two new triangles $\triangle ikl$ and $\triangle jkl$.
This is called a \textbf{surgery by flipping} on the triangulation $\mathcal{T}$,
which is also an isometry of $(S, V)$ with PL metric $w(T)*d_0$.
After the surgery at time $t=T$, we run the $\alpha$-flows on $(S, V, \mathcal{T}')$ with initial metric coming from
the corresponding $\alpha$-flow on $(S, V, \mathcal{T})$ at time $t=T$.

We prove the following result on combinatorial $\alpha$-Yamabe flow and $\alpha$-Calabi flow with surgery.
\begin{theorem}\label{main theorem for alpha Yamabe and Calabi flow with surgery}
Suppose $(S, V)$ is a closed connected marked surface with a PL metric $d_0$ and $\alpha\in \mathbb{R}$ is a constant with
$\alpha\chi(S)\leq 0$.
Then there exists a PL metric in the conformal class $\mathcal{D}(d_0)$ with constant combinatorial $\alpha$-curvature if and only if
one of the following two conditions is satisfied:
\begin{description}
  \item[(1)] The combinatorial $\alpha$-Yamabe flow with surgery
exists for all time and
converges exponentially fast to a PL metric $d^*$ with
constant combinatorial $\alpha$-curvature;
  \item[(2)] The combinatorial $\alpha$-Calabi flow with surgery
exists for all time and
converges exponentially fast to a PL metric $d^*$ with
constant combinatorial $\alpha$-curvature.
\end{description}
\end{theorem}

Applying the discrete maximal principle to the combinatorial $\alpha$-Yamabe flow with surgery suitably, we further have the following result.

\begin{theorem}\label{main theorem convergence of alpha Yamabe and Calabi flow with surgery}
Suppose $(S, V)$ is a closed connected marked surface with a PL metric $d_0$ and $\alpha\in \mathbb{R}$ is a constant with
$\alpha\chi(S)\leq 0$.
Then the combinatorial $\alpha$-Yamabe flow with surgery and the combinatorial $\alpha$-Calabi flow with surgery
exist for all time and
converge exponentially fast to a PL metric $d^*$ with
constant combinatorial $\alpha$-curvature.
\end{theorem}

Theorem \ref{main theorem convergence of alpha Yamabe and Calabi flow with surgery}
confirms another generalized Luo conjecture \cite{L1} on the convergence of the combinatorial Yamabe flow with surgery.
When $\alpha=0$, the convergence of combinatorial Yamabe flow with surgery was proved in \cite{GLSW}
and the convergence of combinatorial Calabi flow with surgery was proved in \cite{ZX}.

The paper is organized as follows.
In Section \ref{section 2}, we prove Theorem \ref{main rigidity intro}
and study the stability of combinatorial $\alpha$-flows on triangulated surfaces.
In Section \ref{section 3}, we prove Theorem \ref{main theorem existence of constant alpha curvature metric},
Theorem \ref{main theorem for alpha Yamabe and Calabi flow with surgery} and
Theorem \ref{main theorem convergence of alpha Yamabe and Calabi flow with surgery}
based on the discrete conformal theory established in \cite{GLSW}.

\section{$\alpha$-curvature and $\alpha$-flows on triangulated surfaces}\label{section 2}
\subsection{Rigidity of $\alpha$-curvature on triangulated surfaces}
Suppose $(S, V, \mathcal{T})$ is a triangulated surface with a PL metric $d$ and $w: V\rightarrow (0, +\infty)$ is a positive function defined on $V$.
Set $h: \mathbb{R}^n_{>0}\rightarrow \mathbb{R}^n$ be the homeomorphism defined by $u_i=h(w_i)=\ln w_i$.
Then $w$ is a conformal factor of $d$ on $(S, V, \mathcal{T})$
if and only if $u: V\rightarrow \mathbb{R}$ is in the following space
\begin{equation}
\mathcal{U}_{ijk}^{\mathcal{T}}(d)\triangleq \{(u_i,u_j,u_k)\in \mathbb{R}^3| \frac{d_{rs}}{e^{u_t}}+\frac{d_{rt}}{e^{u_s}}>\frac{d_{st}}{e^{u_r}}, \{r,s,t\}=\{i,j,k\}\}
\end{equation}
for every triangle $\triangle ijk\in F$.
It is observed by Luo \cite{L1} that the non-convex simply connected space $\mathcal{U}_{ijk}^{\mathcal{T}}(d)$
is the image of the convex space $\{(w_i^{-1}, w_j^{-1}, w_k^{-1})\in \mathbb{R}^3_{>0}|(\frac{d_{rs}}{w_t}, \frac{d_{rt}}{w_s}, \frac{d_{st}}{w_r})\in \Delta\}$
under the homeomorphism $-h$, where $\Delta= \{(x_1, x_2, x_3)\in \mathbb{R}^3_{>0}|x_i+x_j>x_k, \text{where $i,j, k$ are distinct}\}$.
For a nondegenerate triangle $\triangle ijk$ in $(S, V, \mathcal{T})$ with a PL metric $d$, we denote the inner angle at $i$ as $\theta_i$ for simplicity.
Luo \cite{L1} proved the following lemma for a triangle.
\begin{lemma}\label{L_ijk}
The $3\times 3$ matrix $[\frac{\partial \theta_r}{\partial u_s}]_{3\times 3}$ is symmetric, negative semi-definite with null space
$\{(t, t, t)\in \mathbb{R}^3| t\in \mathbb{R}\}$.
\end{lemma}
Lemma \ref{L_ijk} and the simply connectness of $\mathcal{U}_{ijk}^{\mathcal{T}}(d)$ implies that
\begin{equation}\label{Luo's energy function}
F_{ijk}(u)=\int_{u_0}^u \theta_idu_i+\theta_jdu_j+\theta_kdu_k
\end{equation}
is well-defined on $\mathcal{U}_{ijk}^{\mathcal{T}}(d)$. Furthermore, $F_{ijk}(u)$ is locally concave on $\mathcal{U}_{ijk}^{\mathcal{T}}(d)$ and
locally strictly concave on $\mathcal{U}_{ijk}^{\mathcal{T}}(d)\cap \{u_i+u_j+u_k=c\}$.
Bobenko-Pinkall-Springborn \cite{BPS} obtained the explicit formula of $F_{ijk}$ using Milnor's Lobachevsky function and extended
$F_{ijk}$ to be a globally concave function $\widetilde{F}_{ijk}$ defined on $\mathbb{R}^3$.
See also \cite{X4}.
Luo \cite{L2} studied Bobenko-Pinkall-Springborn's extension and obtained a general extension method for similar problems
without involving Milnor's Lobachevsky function, which has lots of applications (see \cite{L2,LY,X1,X2} for example). Here we take Luo's approach.

\begin{lemma}[\cite{L2}]\label{extension of inner angle}
Let $l_1, l_2, l_3$ and $\theta_1, \theta_2, \theta_3$ be the edge lengths and inner angles of a triangle $\triangle$ in $\mathbb{E}^2$,
or $\mathbb{H}^2$, or $\mathbb{S}^2$ so that the $l_i$-th edge is opposite to the angle $\theta_i$.
Consider $\theta_i=\theta_i(l)$ as a function of $l=(l_1, l_2, l_3)$.
\begin{description}
  \item[(1)] If $\triangle$ is Euclidean or hyperbolic, the angle function $\theta_i$ defined on
  $$\Omega=\{(l_1, l_2, l_3)\in \mathbb{R}^3|l_1+l_2>l_3, l_1+l_3>l_2, l_2+l_3>l_1\}$$
  can be extended continuously by constant functions to a function $\widetilde{\theta}_i$ on $\mathbb{R}^3_{>0}$.
  \item[(2)] If $\triangle$ is spherical, the angle function $\theta_i$ defined on
  $$\Omega=\{(l_1, l_2, l_3)\in \mathbb{R}^3|l_1+l_2>l_3, l_1+l_3>l_2, l_2+l_3>l_1, l_1+l_2+l_3<2\pi\}$$
  can be extended continuously by constant functions to a function $\widetilde{\theta}_i$ on $(0, \pi)^3$.
\end{description}
\end{lemma}

Before going on, we recall the following result of Luo in \cite{L2}.
\begin{definition}
A differential 1-form $w=\sum_{i=1}^n a_i(x)dx^i$ in an open set $U\subset \mathbb{R}^n$ is said to be continuous if
each $a_i(x)$ is continuous on $U$.  A continuous differential 1-form $w$ is called closed if $\int_{\partial \tau}w=0$ for each
triangle $\tau\subset U$.
\end{definition}

\begin{theorem}[\cite{L2}, Corollary 2.6]\label{Luo's convex extention}
Suppose $X\subset \mathbb{R}^n$ is an open convex set and $A\subset X$ is an open subset of $X$ bounded by a $C^1$
smooth codimension-1 submanifold in $X$. If $w=\sum_{i=1}^na_i(x)dx_i$ is a continuous closed 1-form on $A$ so that
$F(x)=\int_a^x w$ is locally convex on $A$ and each $a_i$ can be extended continuous to $X$ by constant functions to a
function $\widetilde{a}_i$ on $X$, then  $\widetilde{F}(x)=\int_a^x\sum_{i=1}^n\widetilde{a}_i(x)dx_i$ is a $C^1$-smooth
convex function on $X$ extending $F$.
\end{theorem}

Using Lemma \ref{extension of inner angle} and Theorem \ref{Luo's convex extention}, we have
\begin{lemma}[\cite{BPS, L2}]\label{extension lemma}
The function $F_{ijk}(u)$ in (\ref{Luo's energy function})
could be extended to be a $C^1$-smooth concave function
\begin{equation}\label{extension of F_{ijk}}
\widetilde{F}_{ijk}(u)=\int_{u_0}^u \widetilde{\theta}_idu_i+\widetilde{\theta}_jdu_j+\widetilde{\theta}_kdu_k
\end{equation}
defined for $u\in \mathbb{R}^3$,
where the extension $\widetilde{\theta}_i$ of $\theta_i$ by constant is defined to be
$\widetilde{\theta}_i=\pi$ when $l_{jk}\geq l_{ik}+l_{ij}$ and $\widetilde{\theta}_i=0$
when $l_{ik}\geq l_{jk}+l_{ij}$ or $l_{ij}\geq l_{jk}+l_{ik}$.
\end{lemma}
\textbf{Proof of Theorem \ref{main rigidity intro}:}
The proof is parallelling to that of Theorem 3.3 in \cite{GX5}. For completeness, we give the proof here.
Suppose $w_0\in \Omega^\mathcal{T}(d)$ is a conformal factor and $u_0=\ln w_{0}$.
Then we can define the following Ricci energy  $F(u)$ by $\overline{R}$
\begin{equation}\label{F}
\begin{aligned}
F(u)=-\sum_{\Delta ijk\in F}F_{ijk}+\int_{u_0}^u\sum_{i=1}^N(2\pi-\overline{R}_iw_i^\alpha)du_i.
\end{aligned}
\end{equation}
Note that the function $F_{ijk}$ is smooth on
$\mathcal{U}^{\mathcal{T}}(d)= h(\Omega^{\mathcal{T}}(d))$.
By direct calculations, we have
\begin{equation*}
\operatorname{Hess}_u F= L-\alpha\left(
              \begin{array}{ccc}
                \overline{R}_1w_i^{\alpha} &   &   \\
                  & \ddots &   \\
                  &   & \overline{R}_Nw_N^{\alpha} \\
              \end{array}
            \right),
\end{equation*}
where
\begin{equation}\label{definition of L}
\begin{aligned}
L=(L_{ij})_{N\times N}=\frac{\partial (K_1, \cdots, K_N)}{\partial (u_1, \cdots, u_N)}
=\left(
                                                                                         \begin{array}{ccc}
                                                                                           \frac{\partial K_1}{\partial u_1} & \cdots & \frac{\partial K_1}{\partial u_N} \\
                                                                                           \vdots & \ddots & \vdots \\
                                                                                           \frac{\partial K_N}{\partial u_1} & \cdots & \frac{\partial K_N}{\partial u_N} \\
                                                                                         \end{array}
                                                                                       \right).
\end{aligned}
\end{equation}
The matrix $L$ has the following property \cite{L1}.
\begin{lemma}[\cite{L1}]\label{property of Euclidean L}
For a triangulated surface $(S, V, \mathcal{T})$ with a PL metric $d$,
the matrix $L$ is symmetric and positive semi-definite on $\mathcal{U}^{\mathcal{T}}(d)$ with kernel $\{t \textbf{1}|t\in \mathbb{R}\}$,
where $\textbf{1}=(1, \cdots, 1)$.
\end{lemma}
If $\alpha\overline{R}\equiv 0$,
then $\operatorname{Hess}_u F$ is positive semi-definite with kernel $\{(t, \cdots, t)|t\in \mathbb{R}\}$ and
$F$ is locally convex.
If $\alpha\overline{R}\leq 0$ and $\alpha\overline{R}\not\equiv 0$, then $\operatorname{Hess}_u F$ is positive definite
and $F$ is locally strictly convex.

By Lemma \ref{extension lemma}, $F_{ijk}$ defined on $\mathcal{U}_{ijk}^{\mathcal{T}}(d)$
could be extended to be $\widetilde{F}_{ijk}$ defined by (\ref{extension of F_{ijk}}) on $\mathbb{R}^3\hookrightarrow \mathbb{R}^N$.
And the second term $\int_{u_0}^u\sum_{i=1}^N(2\pi-\overline{R}_iw_i^\alpha)du_i$ in (\ref{F})
can be naturally defined on $\mathbb{R}^N$, then we have the following extension $\widetilde{F}(u)$ defined on $\mathbb{R}^N$
of the Ricci energy function $F(u)$
$$\widetilde{F}(u)=-\sum_{\triangle ijk\in F}\widetilde{F}_{ijk}+\int_{u_0}^u\sum_{i=1}^N(2\pi-\overline{R}_iw_i^\alpha)du_i.$$
As $\widetilde{F}_{ijk}$ is $C^1$-smooth concave by Lemma \ref{extension lemma} and $\int_{u_0}^u\sum_{i=1}^N(2\pi-\overline{R}_iw_i^\alpha)du_i$
is a well-defined convex function on $\mathbb{R}^N$ for $\alpha\overline{R}\leq 0$,
we have $\widetilde{F}(u)$ is a $C^1$-smooth convex function on $\mathbb{R}^N$. Furthermore,
\begin{equation*}
\nabla_{u_i}\widetilde{F}=-\sum_{\triangle ijk\in F}\widetilde{\theta}_i+2\pi-\overline{R}_iw_i^\alpha=\widetilde{K}_i-\overline{R}_iw_i^\alpha,
\end{equation*}
where $\widetilde{K}_i=2\pi-\sum_{\triangle ijk\in F}\widetilde{\theta}_i$.
Then we have $\widetilde{F}(u)$ is convex on $\mathbb{R}^N$ and locally strictly convex on $\mathcal{U}^{\mathcal{T}}(d) \cap\{\sum_{i=1}^Nu_i=0\}$ for $\alpha\overline{R}\equiv0$.
Similarly, $\widetilde{F}(u)$ is convex on $\mathbb{R}^N$ and locally strictly convex on
$\mathcal{U}^{\mathcal{T}}(d)$ for $\alpha\overline{R}\leq0$ and $\alpha\overline{R}\not\equiv0$.

If there are two different conformal factors $\overline{w}_{A}, \overline{w}_{B}$ with the same combinatorial
$\alpha$-curvature $\overline{R}$, then
$\overline{u}_A=\ln \overline{w}_{A}\in \mathcal{U}^{\mathcal{T}}(d)$, $\overline{u}_B=\ln \overline{w}_{B}\in \mathcal{U}^{\mathcal{T}}(d)$
are both critical points of the extended Ricci potential $\widetilde{F}(u)$.
It follows that
$$\nabla \widetilde{F}(\overline{u}_A)=\nabla \widetilde{F}(\overline{u}_B)=0.$$
Set
\begin{equation*}
\begin{aligned}
f(t)=&\widetilde{F}((1-t)\overline{u}_A+t\overline{u}_B)\\
=&\sum_{\triangle ijk\in F}f_{ijk}(t)+\int_{u_0}^{(1-t)\overline{u}_A+t\overline{u}_B}\sum_{i=1}^N(2\pi-\overline{R}_iw_i^\alpha)du_i,
\end{aligned}
\end{equation*}
where
$$f_{ijk}(t)=-\widetilde{F}_{ijk}((1-t)\overline{u}_A+t\overline{u}_B).$$
Then $f(t)$ is a $C^1$-smooth convex function on $[0, 1]$ and $f'(0)=f'(1)=0$, which implies that $f'(t)\equiv 0$ for $t\in [0, 1]$.
Note that $\overline{u}_A$ is in the open set $\mathcal{U}^{\mathcal{T}}(d)$,
there exists $\epsilon>0$ such that $(1-t)\overline{u}_A+t\overline{u}_B\in \mathcal{U}^{\mathcal{T}}(d)$ for $t\in [0, \epsilon]$.
Then $f(t)$ is smooth on $[0, \epsilon]$.

In the case of $\alpha\overline{R}\leq 0$ and $\alpha\overline{R}\not\equiv0$, the strict convexity of $\widetilde{F}(u)$
on $\mathcal{U}^{\mathcal{T}}(d)$ implies that $f(t)$ is strictly convex on $[0, \epsilon]$ and $f'(t)$ is a strictly increasing function on $[0, \epsilon]$.
Then $f'(0)=0$ implies $f'(\epsilon)>0$, which contradicts $f'(t)\equiv 0$ on $[0, 1]$. So there exists at most one conformal factor with combinatorial $\alpha$-curvature $\overline{R}$.

For the case of $\alpha\overline{R}\equiv0$, we have $f(t)$ is $C^1$ convex on $[0, 1]$ and smooth on $[0, \epsilon]$.
$f'(t)\equiv 0$ on $[0, 1]$ implies that $f''(t)\equiv 0$ on $[0, \epsilon]$.
Note that, for $t\in [0, \epsilon]$,
\begin{equation*}
\begin{aligned}
f''(t)=(\overline{u}_A-\overline{u}_B) L  (\overline{u}_A-\overline{u}_B)^T.
\end{aligned}
\end{equation*}
By Lemma \ref{property of Euclidean L}, we have $\overline{u}_A-\overline{u}_B=c\textbf{1}$ for some constant $c\in \mathbb{R}$, which
implies that $\overline{w}_A=e^{c}\overline{w}_B$. So there exists at most one conformal factor
with combinatorial $\alpha$-curvature $\overline{R}$ up to scaling.
\qed

Theorem \ref{main rigidity intro} has a direct corollary.
\begin{corollary}
Suppose $(S, V, \mathcal{T})$ is a triangulated closed surface with a PL metric $d$
and $\alpha\in \mathbb{R}$ is a constant such that $\alpha\chi(S)\leq0$.
Then there exists at most one $u^*\in \mathcal{U}^{\mathcal{T}}(d)$ such that the PL metric $e^{u^*}*d$
has constant combinatorial $\alpha$-curvature (up to scaling for $\alpha\chi(S)=0$).
\end{corollary}


\subsection{Combinatorial Yamabe flow of $\alpha$-curvature on triangulated surfaces}

%

By direct calculations, we have the following properties of combinatorial $\alpha$-Yamabe flow.

\begin{lemma}\label{basic property of alpha Yamabe flow}
If $\alpha=0$, $\sum_{i=1}^N u_i$ is invariant along the normalized combinatorial $\alpha$-Yamabe flow (\ref{alpha Yamabe flow}).
If $\alpha\neq 0$, $||w||_\alpha^\alpha=\sum_{i=1}^Nw_i^\alpha$ is invariant
along the normalized combinatorial $\alpha$-Yamabe flow (\ref{alpha Yamabe flow}).
\end{lemma}
For simplicity, we denote the hypersurface invariant along the normalized combinatorial $\alpha$-Yamabe flow (\ref{alpha Yamabe flow}) in Lemma
\ref{basic property of alpha Yamabe flow} as $P$ in the following. Note that the hypersurface $P$ is determined by the initial value $u(0)$.

\begin{theorem}\label{stability of alpha Yamabe flow}
Suppose $d_0$ is a PL metric on a triangulated surface $(S, V, \mathcal{T})$ and $\alpha\in \mathbb{R}$.
If the solution of normlized combinatorial $\alpha$-Yamabe flow (\ref{alpha Yamabe flow})
on $(S, V, \mathcal{T})$
converges, then the limit metric is a constant combinatorial $\alpha$-curvature PL metric.
Furthermore, suppose there exists a constant combinatorial $\alpha$-curvature PL metric $d^*=e^{u^*}*d_0$
on a triangulated surface $(S, V, \mathcal{T})$ with $\alpha\chi(S)\leq 0$,
there exists a constant $\delta>0$
such that if  $||R_\alpha(u(0))-R_\alpha(u^*)||<\delta$  and $u^*\in P$,
then the combinatorial $\alpha$-Yamabe flow (\ref{alpha Yamabe flow})
on $(S, V, \mathcal{T})$ exists for all time and converges exponentially fast to $u^*$.
\end{theorem}

\proof
Suppose $u(t)$ is a solution of the normalized combinatorial
$\alpha$-Yamabe flow (\ref{alpha Yamabe flow}).
If $u(\infty)=\lim_{t\rightarrow +\infty}u(t)$ exists in $\mathcal{U}^\mathcal{T}(d)$, then we have
$R_{\alpha}(u(\infty))=\lim_{t\rightarrow +\infty}R_\alpha(u(t))$ exists.
Furthermore, there exists $\xi_n\in (n, n+1)$
such that
$$u_i(n+1)-u_i(n)=u'_i(\xi_n)=R_{\alpha, av}-R_{\alpha, i}(\xi_n)\rightarrow 0,$$
which implies that $R_{\alpha}(u(\infty))=R_{av}$ and $u(\infty)*d_0$ is a constant $\alpha$-curvature PL metric.

Suppose $u^*$ corresponds to a constant $\alpha$-curvature metric.
Set $\Gamma_i(u)=R_{\alpha, av}-R_{\alpha, i}$. By direct calculations, we have
\begin{equation*}
\begin{aligned}
\frac{\partial \Gamma_i}{\partial u_j}|_{u=u^*}
=&-\frac{1}{w_i^\alpha}\frac{\partial K_i}{\partial u_j}+\alpha R_{\alpha, av}(\delta_{ij}-\frac{w_j^\alpha}{||w||_\alpha^\alpha})\\
=&\alpha R_{\alpha, av}\delta_{ij}-\frac{1}{w_i^\alpha}(\frac{\partial K_i}{\partial u_j}+\alpha R_{\alpha, av}\frac{w_i^\alpha w_j^\alpha}{||w||_\alpha^\alpha}).
\end{aligned}
\end{equation*}
Set $w^\alpha=(w_1^\alpha, \cdots, w_N^\alpha)^T$ and $\Sigma=diag\{w_1, \cdots, w_N\}$, then
\begin{equation}\label{gamma matrix}
\begin{aligned}
(\frac{\partial\Gamma}{\partial u})|_{u=u^*}
=&\alpha R_{\alpha, av}I-\Sigma^{-\alpha}(L+\alpha R_{\alpha, av}\frac{w^\alpha\cdot (w^\alpha)^T}{||w||_\alpha^\alpha})\\
=&-\Sigma^{-\alpha/2}\left(\Lambda_{\alpha}-\alpha R_{\alpha, av}[I-\frac{w^{\alpha/2}\cdot (w^{\alpha/2})^T}{||w||_\alpha^\alpha}]\right)\Sigma^{\alpha/2},
\end{aligned}
\end{equation}
where $\Lambda_{\alpha}=\Sigma^{-\alpha/2}L\Sigma^{-\alpha/2}$.
Note that the matrix $I-\frac{w^{\alpha/2}\cdot (w^{\alpha/2})^T}{||w||_\alpha^\alpha}$ has
eigenvalues 1 ($N-1$ times) and 0 (1 time) and kernel $\{cw^{\frac{\alpha}{2}}| c\in \mathbb{R}\}$ and $\Lambda_\alpha$
is positive semi-definite with $1$-dimensional kernel $\{cw^{\frac{\alpha}{2}}| c\in \mathbb{R}\}$.
Then if the first nonzero eigenvalue $\lambda_1(\Lambda_\alpha)$ of $\Lambda_\alpha$ satisfies
$\lambda_1(\Lambda_\alpha)>\alpha R_{\alpha, av}$, especially if $\alpha R_{\alpha, av}\leq 0$,
we have $(\frac{\partial\Gamma}{\partial u})|_{u=u^*}$ has $N-1$ negative eigenvalues and a zero eigenvalue with eigenspace $\{(c,c,\cdots,c)\in \mathbb{R}^N| c\in \mathbb{R}\}$.

If $\alpha\neq 0$, set $u_i=e^{-\frac{\alpha}{2}u_i^*}\widetilde{u}_i$ and $\Sigma^*=\Sigma|_{u=u^*}$,
then $u=(\Sigma^*)^{-\frac{\alpha}{2}}\widetilde{u}$ in matrix form.
The combinatorial $\alpha$-Yamabe flow (\ref{alpha Yamabe flow})
could be written as
\begin{equation}\label{CYF in tilde u}
  \frac{d\widetilde{u}}{dt}=(\Sigma^*)^{\frac{\alpha}{2}} \Gamma((\Sigma^*)^{-\frac{\alpha}{2}}\widetilde{u})
\end{equation}
in the variable $\widetilde{u}$.
Set $\widetilde{\Gamma}(\widetilde{u})=(\Sigma^*)^{\frac{\alpha}{2}} \Gamma((\Sigma^*)^{-\frac{\alpha}{2}}\widetilde{u})$ and
$\widetilde{u}^*=(\Sigma^*)^{\frac{\alpha}{2}}u^*$. Then
\begin{equation*}
(\frac{\partial\widetilde{\Gamma}}{\partial \widetilde{u}})|_{\widetilde{u}=\widetilde{u}^*}
=(\Sigma^*)^{\frac{\alpha}{2}}\cdot (\frac{\partial\Gamma}{\partial u})|_{u=u^*}\cdot (\Sigma^*)^{-\frac{\alpha}{2}},
\end{equation*}
which is symmetric by (\ref{gamma matrix}) and
 negative semi-definite with kernel
$\{c(e^{\frac{\alpha}{2}u_1^*}, \cdots, e^{\frac{\alpha}{2}u_N^*})|c\in \mathbb{R}\}$ under the condition $\alpha R_{\alpha, av}\leq 0$
 by the property of $(\frac{\partial\Gamma}{\partial u})|_{u=u^*}$.

Note that the hypersurface $P$ invariant along the
combinatorial $\alpha$-Yamabe flow (\ref{alpha Yamabe flow}) could be written as
\begin{equation}\label{invariant hypersurface}
\sum_{i=1}^N e^{\alpha e^{-\frac{\alpha}{2}u_i^*}\widetilde{u}_i}=\sum_{i=1}^N e^{\alpha e^{-\frac{\alpha}{2}u_i^*}\widetilde{u}_i(0)}
\end{equation}
in $\widetilde{u}$.
The normal space of the hypersurface (\ref{invariant hypersurface}) in $\widetilde{u}$ at $\widetilde{u}^*$ is $\{c(e^{\frac{\alpha}{2}u_1^*}, \cdots, e^{\frac{\alpha}{2}u_N^*})|c\in \mathbb{R}\}$, which is exactly the kernel space
of $(\frac{\partial\widetilde{\Gamma}}{\partial \widetilde{u}})|_{\widetilde{u}=\widetilde{u}^*}$.
Therefore, $(\frac{\partial\widetilde{\Gamma}}{\partial \widetilde{u}})|_{\widetilde{u}=\widetilde{u}^*}$
is negative definite on the tangential space of $P$ at $\widetilde{u}^*$ and
$\widetilde{u}^*$ is a local attractor of the combinatorial $\alpha$-Yamabe flow (\ref{CYF in tilde u}) in $\widetilde{u}$.
As a result, the local convergence of the solution $\tilde{u}(t)$ of the combinatorial $\alpha$-Yamabe flow (\ref{CYF in tilde u}) follows from the Lyapunov stability theorem (\cite{P}, Chapter 5), from which the local convergence of the solution of the
combinatorial $\alpha$-Yamabe flow (\ref{alpha Yamabe flow}) follows.

The proof for the case of $\alpha=0$ is similar.
\qed

\subsection{Combinatorial Calabi flow of $\alpha$-curvature on triangulated surfaces}
Similar to the combinatorial $\alpha$-Yamabe flow, we have the following properties of combinatorial $\alpha$-Calabi flow.

\begin{lemma}\label{basic property of alpha Calabi flow}
If $\alpha=0$, $\sum_{i=1}^N u_i$ is invariant along the combinatorial $\alpha$-Calabi flow (\ref{alpha Calabi flow}).
If $\alpha\neq 0$, $||w||_\alpha^\alpha=\sum_{i=1}^Nw_i^\alpha$ is invariant
along the combinatorial $\alpha$-Calabi flow (\ref{alpha Calabi flow}).
\end{lemma}
%

\begin{theorem}\label{stability of alpha Calabi flow}
Suppose $d_0$ is a PL metric on a triangulated surface $(S, V, \mathcal{T})$ and $\alpha\in \mathbb{R}$.
If the solution of combinatorial $\alpha$-Calabi flow
on $(S, V, \mathcal{T})$
converges, then the limit metric is a constant combinatorial $\alpha$-curvature PL metric.
Furthermore, suppose there exists a constant combinatorial $\alpha$-curvature PL metric $d^*=e^{u^*}*d_0$
on $(S, V, \mathcal{T})$ with $\alpha \chi(S)\leq 0$,
there exists a constant $\delta>0$
such that if  $||R_\alpha(u(0))-R_\alpha(u^*)||<\delta$  and $u^*\in P$,
then the combinatorial $\alpha$-Calabi flow (\ref{alpha Calabi flow})
on $(S, V, \mathcal{T})$ exists for all time and converges exponentially fast to $u^*$.
\end{theorem}
\proof
The proof of Theorem \ref{stability of alpha Calabi flow} is similar to that of Theorem \ref{stability of alpha Yamabe flow},
we just give some key calculations.
Set $\Gamma_i(u)=(\Delta_\alpha^\mathcal{T}R_{\alpha})_i$, then
\begin{equation*}
\begin{aligned}
\frac{\partial \Gamma_i}{\partial u_j}|_{u=u^*}
=&-\frac{1}{w_i^\alpha}\sum_{k=1}^NL_{ik}\frac{1}{w_k^\alpha}L_{kj}+\alpha R_{\alpha,av}\frac{1}{w_i^\alpha}L_{ij}.
\end{aligned}
\end{equation*}
In matrix form, we have
\begin{equation*}
\begin{aligned}
(\frac{\partial\Gamma}{\partial u})|_{u=u^*}
=&-\Sigma^{-\alpha}L\Sigma^{-\alpha}L+\alpha R_{\alpha, av}\Sigma^{-\alpha}L\\
=&-\Sigma^{-\alpha/2}\left(\Sigma^{-\alpha/2}L\Sigma^{-\alpha}L\Sigma^{-\alpha/2}-\alpha R_{\alpha, av}\Sigma^{-\alpha/2}L\Sigma^{-\alpha/2}\right)\Sigma^{\alpha/2}.
\end{aligned}
\end{equation*}
If $\alpha \chi(S)\leq 0$, then $(\frac{\partial\Gamma}{\partial u})|_{u=u^*}$ has $N-1$ negative eigenvalue and a zero eigenvalue with
$1$-dimensional kernel $\{(c,c,\cdots,c)\in \mathbb{R}^N| c\in \mathbb{R}\}$.
The following of the proof is the same as that for Theorem \ref{stability of alpha Yamabe flow}, we omit the details here.
\qed

\section{$\alpha$-curvature and $\alpha$-flows on discrete Riemann surfaces}\label{section 3}

Theorem \ref{stability of alpha Yamabe flow} and Theorem \ref{stability of alpha Calabi flow}
gives the long time existence and convergence of the combinatorial $\alpha$-Yamabe flow (\ref{alpha Yamabe flow}) and combinatorial
$\alpha$-Calabi flow (\ref{alpha Calabi flow}) for initial PL metrics with small initial energy.
However, for general initial PL metrics, the combinatorial $\alpha$-Yamabe flow and combinatorial
$\alpha$-Calabi flow may develop singularities, including the conformal factor tends to infinity and some triangle degenerates
along the combinatorial $\alpha$-Yamabe flow and combinatorial $\alpha$-Calabi flow.
To handle the possible singularities along the $\alpha$-flows, we do surgery on the flows by flipping as described in Section \ref{section 1}.


To analyze the behavior of the $\alpha$-flows with surgery, we need to use the discrete conformal theory established
by Gu-Luo-Sun-Wu \cite{GLSW} for PL metrics.
In the following, we briefly recall some results in \cite{GLSW}.
For details of the theory, please refer to Gu-Luo-Sun-Wu's important work \cite{GLSW}.

\subsection{Gu-Luo-Sun-Wu's work on discrete uniformization theorem}
\begin{definition}[\cite{GLSW} Definition 1.1]\label{GLSW's definition of edcf}
Two PL metrics $d, d'$ on $(S, V)$ are discrete conformal if
there exist sequences of PL metrics $d_1=d$, $\cdots$,  $d_m=d'$
on $(S, V)$ and triangulations $\mathcal{T}_1, \cdots, \mathcal{T}_m$ of
$(S, V)$ satisfying
\begin{description}
  \item[(a)] (Delaunay condition) each $\mathcal{T}_i$ is Delaunay in $d_i$,
  \item[(b)] (Vertex scaling condition) if $\mathcal{T}_i=\mathcal{T}_{i+1}$, there exists a function
  $u:V\rightarrow \mathbb{R}$ so that if $e$ is an edge in $\mathcal{T}_i$ with end points $v$ and $v'$,
  then the lengths $l_{d_{i+1}}(e)$ and  $l_{d_{i}}(e)$ of $e$ in $d_i$ and $d_{i+1}$ are related by
  $$l_{d_{i+1}}(e)=l_{d_{i}}(e)e^{u(v)+u(v')},$$
  \item[(c)] if $\mathcal{T}_i\neq \mathcal{T}_{i+1}$, then $(S, d_i)$ is isometric to $(S, d_{i+1})$
  by an isometry homotopic to identity in $(S, V)$.
\end{description}
The discrete conformal class of a PL metric is called a discrete Riemann surface.
\end{definition}
The space of PL metrics on $(S, V)$ discrete conformal to $d$ is called the conformal class of $d$ and
denoted by $\mathcal{D}(d)$.

%

The following discrete uniformization theorem was established in \cite{GLSW}.

\begin{theorem}[\cite{GLSW} Theorem 1.2]\label{Euclidean discrete uniformization}
Suppose $(S, V)$ is a closed connected marked surface and $d$ is a PL metric on $(S, V)$.
Then for any $K^*: V\rightarrow (-\infty, 2\pi)$ with $\sum_{v\in V}K^*(v)=2\pi\chi(S)$,
there exists a PL metric $d'$, unique up to scaling and isometry homotopic to the identity
 on $(S, V)$, such that $d'$ is discrete conformal to $d$ and the discrete curvature
 of $d'$ is $K^*$.
\end{theorem}

Denote the Teichim\"{u}ller space of all PL metrics on $(S, V)$ by $T_{PL}(S, V)$ and decorated
Teichim\"{u}ller space of all equivalence class of decorated hyperbolic metrics on $S-V$ by $T_D(S-V)$.
In the proof of Theorem \ref{Euclidean discrete uniformization}, Gu-Luo-Sun-Wu proved the following result.

\begin{theorem}[\cite{GLSW}]\label{main result of GLSW}
There is a $C^1$-diffeomorphism $\mathbf{A}: T_{PL}(S, V)\rightarrow T_D(S, V)$ between $T_{PL}(S, V)$ and $T_D(S-V)$.
Furthermore, the space $\mathcal{D}(d)\subset T_{PL}(S, V)$ of all equivalence classes of PL metrics
discrete conformal to $d$ is $C^1$-diffeomorphic to $\{p\}\times \mathbb{R}^V_{>0}$ under the diffeomorphism $\mathbf{A}$,
where $p$ is the unique hyperbolic metric on $S-V$ determined by the PL metric $d$ on $(S, V)$.
\end{theorem}

Set $u_i=\ln w_i$ for $w=(w_1, w_2, \cdots, w_n)\in \mathbb{R}^n_{>0}$.
Using the map $\mathbf{A}$, Gu-Luo-Sun-Wu defined the curvature map
\begin{equation}\label{GLSW's curvature map}
\begin{aligned}
\mathbf{F}:\mathbb{R}^n&\rightarrow (-\infty, 2\pi)^n\\
u&\mapsto K_{\mathbf{A}^{-1}(p, w(u))}
\end{aligned}
\end{equation}
and proved the following property of $\mathbf{F}$.

\begin{proposition}[\cite{GLSW}] \label{energy function Euclidean}
\begin{enumerate}
  \item For any $k\in \mathbb{\mathbb{R}}$, $\mathbf{F}(v+k(1, 1, \cdots, 1))=\mathbf{F}(v)$.
  \item There exists a $C^2$-smooth convex function $W: \mathbb{R}^n\rightarrow \mathbb{R}$
so that its gradient $\nabla W$ is $\mathbf{F}$ and the restriction $W: \{u\in \mathbb{R}^n|\sum_{i=1}^nu_i=0\}\rightarrow \mathbb{R}$
is strictly convex.
\end{enumerate}
\end{proposition}

Theorem \ref{main result of GLSW} implies that the union of the admissible spaces $\Omega^{\mathcal{T}}_{D}(d')$
of conformal factors such that $\mathcal{T}$ is Delaunay for $d'\in \mathcal{D}(d)$
is $\mathbb{R}^n_{>0}$. Furthermore, $\mathbf{F}$, which is defined on $\mathbb{R}^n_{>0}$, is
a $C^1$-extension of the curvature $K$ defined on
the space of conformal factors $\Omega^{\mathcal{T}}_{D}(d')$ for $d'\in \mathcal{D}(d)$.
Then we can extend the Euclidean discrete $\alpha$-Laplace operator to be defined on $\mathbb{R}^n_{>0}$, which is the space of
the conformal factors for the discrete conformal class $\mathcal{D}(d)$.
\begin{definition}\label{Euclid Delaunay Laplacian}
Suppose $(S, V)$  is a marked surface with a PL metric $d_0$,
For a function $f: V\rightarrow \mathbb{R}$ on the vertices,
the discrete conformal $\alpha$-Laplace operator of $d\in \mathcal{D}(d_0)$ on $(S, V)$ is
defined to be the map
\begin{equation*}
\begin{aligned}
\Delta_\alpha: \mathbb{R}^V&\longrightarrow \mathbb{R}^V\\
f&\mapsto \Delta_\alpha f,
\end{aligned}
\end{equation*}
where the value of $\Delta_\alpha f$ at $v_i$ is
\begin{equation}\label{alpha Laplacian}
\begin{aligned}
\Delta_\alpha f_i=\frac{1}{w_i^\alpha}\sum_{j\sim i}(-\frac{\partial \mathbf{F}_i}{\partial u_j})(f_j-f_i)=-\frac{1}{w_i^\alpha}(\widetilde{L}f)_i,
\end{aligned}
\end{equation}
where $(p, w)=\mathbf{A}(d)$ and $\widetilde{L}_{ij}=\frac{\partial \mathbf{F}_i}{\partial u_j}$ is
an extension of $L_{ij}=\frac{\partial K_i}{\partial u_j}$ for $u=\ln w\in\mathcal{U}^{\mathcal{T}}_{D}(d')=\ln \Omega^{\mathcal{T}}_{D}(d')$, $d'\in \mathcal{D}(d)$.
\end{definition}

\begin{remark}
Note that $\mathbf{F}$ is $C^1$-smooth in $u\in\mathbb{R}^n$ and $\Delta^{\mathcal{T}}$ is independent of
the Delaunay triangulations of a PL metric, so
the operator $\Delta_\alpha$ is well-defined on $\mathbb{R}^n$.
Furthermore, $\Delta_\alpha$ is continuous and piecewise smooth on $\mathbb{R}^n$ as a matrix-valued function of $u$ (\cite{GLSW}, Lemma 5.1).
\end{remark}
\subsection{Rigidity of $\alpha$-curvature on discrete Riemannian surfaces}
Following Gu-Luo-Sun-Wu's approach, we can define the $\alpha$-curvature  on discrete Riemannian surfaces as follows.

\begin{definition}\label{def of alpha curv on dis Riem surf}
Suppose $(S, V)$ is a marked closed surface with a PL metric $d$, $\alpha\in \mathbb{R}$ is a constant and
$\mathbf{F}$ is the curvature map in (\ref{GLSW's curvature map}).
The $\alpha$-curvature on the discrete Riemann surface $\mathcal{D}(d)$ is defined to be
\begin{equation}
\mathbf{F}_{\alpha, i}=\frac{\mathbf{F}_i}{w_i^\alpha}.
\end{equation}
\end{definition}

\begin{remark}
Note that $\alpha$-curvature on a discrete Riemann surface is well-defined and
an extension of the combinatorial $\alpha$-curvature on a triangulated surface.
If $\mathcal{T}$ is a Delaunay triangulation of the marked surface $(S, V)$,
$\Omega^{\mathcal{T}}_{D}(d')$ is the space
of conformal factors such that $\mathcal{T}$ is Delaunay for $d'\in \mathcal{D}(d)$, then
$\mathbf{F}_\alpha|_{\mathcal{U}^{\mathcal{T}}_{D}(d')}=R_\alpha$.
\end{remark}

Denote the space of conformal factors by $\Omega(d)$ and set $\mathcal{U}(d)=\ln\Omega(d)$.
Similar to Theorem \ref{main rigidity intro} for $\alpha$-curvature on triangulated surfaces,
we have the following global rigidity for $\alpha$-curvature on discrete Riemann surfaces.

\begin{theorem}\label{rigidity of alpha curvature map}
Suppose $(S, V)$ is a marked surface with a PL metric $d$ and $\alpha\in \mathbb{R}$ is a constant with $\alpha\chi(S)\leq 0$.
$\overline{F}$ is a function defined on the vertices.
\begin{description}
  \item[(1)] If $\alpha \overline{F}\equiv 0$, then there exists at most one conformal factor $u^*\in\mathcal{U}(d)$ up to scaling such that
 $\mathbf{A}^{-1}(p, w(u^*))\in \mathcal{D}(d)$ has combinatorial $\alpha$-curvature $\overline{F}$.
  \item[(2)] If $\alpha \overline{F}\leq 0$ and $\alpha \overline{F}\not\equiv 0$, then there exists at most one conformal factor $u^*\in\mathcal{U}(d)$ such that $\mathbf{A}^{-1}(p, w(u^*))\in \mathcal{D}(d)$ has combinatorial $\alpha$-curvature $\overline{F}$.
\end{description}
\end{theorem}
\proof
Define the energy function
\begin{equation}
W_\alpha(u)=W(u)-\int_0^u\sum_{i=1}^N \overline{F}_iw_i^\alpha du_i.
\end{equation}
By Proposition \ref{energy function Euclidean}, $W_\alpha$ is a well-defined $C^2$-smooth function defined on $\mathbb{R}^n$.
Furthermore, we have
$$\nabla_{u_i} W_\alpha=\mathbf{F}_i-\overline{F}_iw_i^\alpha.$$
$\mathbf{A}^{-1}(p, w(u^*))\in \mathcal{D}(d)$ for $u^*\in \mathcal{U}(d)$ has combinatorial $\alpha$-curvature $\overline{F}$ if and only if $\nabla W_\alpha(u^*)=0$.
By direct calculations, we have
\begin{equation*}
\begin{aligned}
\operatorname{Hess} W_\alpha=L-\alpha\left(
                          \begin{array}{ccc}
                            \overline{F}_1w_1^\alpha &   &   \\
                              & \ddots &   \\
                              &   & \overline{F}_Nw_N^\alpha \\
                          \end{array}
                        \right).
\end{aligned}
\end{equation*}
If $\alpha \overline{F}\equiv0$, $\operatorname{Hess} W_\alpha$ is positive semi-definite with kernel $t\mathbf{1}=(t,\cdots, t)$
and $W_\alpha|_{\Sigma_0}$ is a strictly convex function on $\Sigma_0=\{u_1+\cdots+u_N=0\}$.
If $\alpha \overline{F}\leq 0$ and $\alpha \overline{F}\not\equiv 0$, then $\operatorname{Hess} W_\alpha$ is positive definite
and $W_\alpha$
is strictly convex on $\mathbb{R}^n$.

Recall the following well known fact from analysis.
\begin{lemma}\label{gradient of convex function}
If $W: \Omega\rightarrow \mathbb{R}$ is a $C^1$-smooth strictly convex function on an open convex set $\Omega\subset \mathbb{R}^m$,
then its gradient $\nabla W: \Omega\rightarrow \mathbb{R}^m$ is an embedding.
\end{lemma}
Then the rigidity follows from Lemma \ref{gradient of convex function}.
\qed

\begin{corollary}
Suppose $(S, V)$ is a marked surface with a PL metric $d$ and $\alpha\in \mathbb{R}$ is a constant with $\alpha\chi(S)\leq 0$.
Then the constant combinatorial $\alpha$-curvature PL metric in $\mathcal{D}(d)$ is unique
(up to scaling if $\alpha\chi(S)=0$).
\end{corollary}

\subsection{Combinatorial $\alpha$-Yamabe flow with surgery}
By Gu-Luo-Sun-Wu's discrete conformal theory \cite{GLSW},
the normalized combinatorial $\alpha$-Yamabe flow with surgery takes the following form.

\begin{definition}
Suppose $(S, V)$ is a marked surface with a PL metric $d_0$. The combinatorial $\alpha$-Yamabe flow with surgery is defined to be
\begin{equation}\label{alpha Yamabe flow with surgery}
\begin{aligned}
\left\{
  \begin{array}{ll}
    \frac{du_i}{dt}=\mathbf{F}_{\alpha, av}-\mathbf{F}_{\alpha,i}, & \hbox{ } \\
    u_i(0)=0, & \hbox{ }
  \end{array}
\right.
\end{aligned}
\end{equation}
where
$\mathbf{F}_{\alpha, av}=\frac{2\pi\chi(S)}{\sum_{i=1}^Nw_i^\alpha}.$
\end{definition}

It is straightforward to check that $\sum_{i=1}^N w_i^\alpha$ ($\sum_{i=1}^Nu_i$ for $\alpha=0$) is invariant along the
combinatorial $\alpha$-Yamabe flow with surgery (\ref{alpha Yamabe flow with surgery}).

Similar to the results in \cite{CL,Ge, GX4}, we have the following result for combinatorial $\alpha$-Yamabe flow
with surgery.

\begin{theorem}\label{equivalene for alpha Yamabe flow}
Suppose $(S, V)$ is a closed connected marked surface with a PL metric $d_0$.
$\alpha\in \mathbb{R}$ is a constant such that $\alpha\chi(S)\leq 0$.
Then there exists a constant $\alpha$-curvature PL metric in $\mathcal{D}(d_0)$
if and only if the combinatorial $\alpha$-Yamabe flow with surgery (\ref{alpha Yamabe flow with surgery}) exists for all time
and converges to some $u^*\in \mathcal{U}(d_0)$.
\end{theorem}
\proof
When $\alpha=0$, the combinatorial $\alpha$-Yamabe flow with surgery (\ref{alpha Yamabe flow with surgery}) is the
Yamabe flow with surgery studied in \cite{L1, GLSW}, where the conclusion has been proved.
We only prove the case $\alpha\neq 0$ here.

If the solution $u(t)$ of combinatorial $\alpha$-Yamabe flow with surgery (\ref{alpha Yamabe flow with surgery})
converges to $u^*\in \mathcal{U}(d_0)$, then
we have $\mathbf{F}_{\alpha}(u^*)=\lim_{t\rightarrow +\infty}\mathbf{F}_{\alpha}(u(t))$
by the $C^1$ smoothness of $\mathbf{F}$. For any $n\in \mathbb{N}$, there exists $\xi_n\in (n, n+1)$ such
that
$$u_i(n+1)-u_i(n)=u_i'(\xi_n)=\mathbf{F}_{\alpha, av}-\mathbf{F}_{\alpha, i}(u(\xi_n)).$$
Set $n\rightarrow +\infty$, then we have
$$\mathbf{F}_{\alpha, i}(u^*)=\lim_{n\rightarrow +\infty}\mathbf{F}_{\alpha, i}(u(\xi_n))=\mathbf{F}_{\alpha, av},$$
which implies that $u^*$ is a conformal factor in $\mathcal{U}(d_0)$ with constant combinatorial $\alpha$-curvature.

Conversely, suppose $u^*$ is a conformal factor in $\mathcal{U}(d_0)$ with constant combinatorial $\alpha$-curvature.
Then the constant curvature must be the constant $\mathbf{F}_{\alpha, av}=\frac{2\pi\chi(S)}{\sum_{i=1}^N w_i^\alpha}$.
Set
$$W_\alpha(u)=W(u)-\mathbf{F}_{\alpha, av}\int_{u^*}^{u}\sum_{i=1}^Nw_i^\alpha du_i.$$
Then $W_\alpha$ is a well-defined $C^2$-smooth convex function defined on $\mathbb{R}^n$ and $W_\alpha(u+k\mathbf{1})=W_\alpha(u)$.
Note that $\nabla W_\alpha(u^*)=0$, we have $\lim_{u\rightarrow \infty}W_\alpha(u)|_{P}=+\infty$, where
$P=\{u\in \mathbb{R}^n|\sum_{i=1}^N w_i^\alpha=N\}$. This implies that $W_\alpha(u)|_{P}$
is a proper function on $P$.

Note that
\begin{equation*}
\begin{aligned}
\frac{d(W_{\alpha}(u(t)))}{dt}
=&\sum_{i=1}^N\frac{\partial W_\alpha}{\partial u_i}\cdot \frac{du_i}{dt}
=\sum_{i=1}^N(\mathbf{F}_i-\mathbf{F}_{\alpha, av}w_i^\alpha)(\mathbf{F}_{\alpha, av}-\mathbf{F}_{\alpha, i})\\
=&-\sum_{i=1}^N(\mathbf{F}_{\alpha, av}-\mathbf{F}_{\alpha, i})^2w_i^\alpha\leq 0.
\end{aligned}
\end{equation*}
So we have
$0\leq W_{\alpha}(u(t))\leq W_{\alpha}(u(0))$.
Note that $\sum_{i=1}^N w_i^\alpha$ is invariant along
the combinatorial $\alpha$-Yamabe flow with surgery (\ref{alpha Yamabe flow with surgery}), we have
the solution $u(t)$ of the combinatorial $\alpha$-Yamabe flow with surgery lies in a compact subset of $P$ by the properness of
$W_\alpha$ on $P$. Then the solution of the combinatorial $\alpha$-Yamabe flow with surgery (\ref{alpha Yamabe flow with surgery}) exists for all time
and $\lim_{t\rightarrow +\infty}W_{\alpha}(u(t))$ exists. Furthermore,
\begin{equation*}
\begin{aligned}
0=&\lim_{n\rightarrow +\infty}(W_\alpha(u(n+1)-W_\alpha(u(n))))
=\lim_{n\rightarrow +\infty}\frac{dW_\alpha(u(t))}{dt}|_{t=\xi_n}\\
=&-\lim_{n\rightarrow +\infty}\sum_{i=1}^N(\mathbf{F}_{\alpha, av}-\mathbf{F}_{\alpha, i})^2w_i^\alpha|_{t=\xi_n}.
\end{aligned}
\end{equation*}
Then we have $\lim_{n\rightarrow +\infty} \mathbf{F}_{\alpha}(u(\xi_n))=\mathbf{F}_{\alpha, av}=\mathbf{F}_{\alpha}(u^*)$,
which implies that $\lim_{n\rightarrow +\infty} u(\xi_n)=u^*$ by Theorem \ref{rigidity of alpha curvature map}.

Set $\Gamma_i(u)=\mathbf{F}_{\alpha, av}-\mathbf{F}_{\alpha,i}$. Similar to the proof of Theorem \ref{stability of alpha Yamabe flow},
we can check that $(\frac{\partial\Gamma}{\partial u})|_{u^*}$ has $N-1$ negative eigenvalue and a zero eigenvalue.
The local convergence of the solution of combinatorial $\alpha$-Yamabe flow with surgery to $u^*$ follows
by the same argument in the proof of Theorem \ref{stability of alpha Yamabe flow}.
The conclusion then follows by the local convergence and  $\lim_{n\rightarrow +\infty} u(\xi_n)=u^*$.
\qed

\begin{remark}
The proof of Theorem \ref{equivalene for alpha Yamabe flow} suggests the following generalization of
the combinatorial $\alpha$-Yamabe flow with surgery.
Suppose $\overline{F}$ is a function defined on the vertices, then $\overline{F}$ is the combinatorial
$\alpha$-curvature of a PL metric in $\mathcal{D}(d_0)$ if and only if the solution of combinatorial
$\alpha$-Yamabe flow with surgery for $\overline{F}$ (defined similarly) exists for all time and converges to a conformal
factor $u^*\in \mathcal{U}(d_0)$.
\end{remark}

In the case of $\alpha\chi(S)\leq0$, we can further prove the existence of constant combinatorial $\alpha$-curvature metric and then
obtain a generalization of the discrete uniformization theorem obtained in \cite{GLSW}.

By direct calculations,
the curvature $\mathbf{F}_{\alpha,i}$
evolves according to the following equation
\begin{equation}\label{evolution of curv under alpha Yamabe flow}
\begin{aligned}
\frac{d\mathbf{F}_{\alpha, i}}{dt}=(\Delta_\alpha\mathbf{F}_{\alpha})_i+\alpha \mathbf{F}_{\alpha,i}(\mathbf{F}_{\alpha,i}-\mathbf{F}_{\alpha,av})
\end{aligned}
\end{equation}
along the combinatorial $\alpha$-Ricci flow.

Note that the surgery ensures that the weight
$$\omega_{ij}=\frac{1}{w_i^\alpha}\frac{\partial \mathbf{F}_i}{\partial u_j}=\frac{\cot\theta_{k}^{ij}+\cot\theta_{l}^{ij}}{w_i^\alpha}\geq 0$$
along the combinatorial $\alpha$-Yamabe flow with surgery (\ref{alpha Yamabe flow with surgery}).
This motives us to use the following discrete maximal principle. The readers can refer to \cite{GX4} for a proof.
\begin{theorem}\label{Maximum priciple}(Maximum Principle)
Let $f: V\times [0, T)\rightarrow \mathbb{R}$ be a $C^1$ function such that
$$\frac{\partial f_i}{\partial t}\geq \Delta f_i+ \Phi_i(f_i), \ \ \forall (i, t)\in V\times [0, T), $$
where the Laplacian operator is defined as
$$\Delta f_i=\sum_{j\sim i}a_{ij}(t)(f_j-f_i)$$
with $a_{ij}\geq 0$
and
$\Phi_i: \mathbb{R}\rightarrow \mathbb{R}$ is a local Lipschitz function.
Suppose there exists $C_1\in \mathbb{R}$ such that $f_i(0)\geq C_1$ for all $i\in V$. Let $\varphi$ be the
solution to the associated ODE
\begin{equation*}
\begin{aligned}
\left\{
  \begin{array}{ll}
    \frac{d \varphi}{dt}=\Phi_i(\varphi), \\
    \varphi(0)=C_1,
  \end{array}
\right.
\end{aligned}
\end{equation*}
then
$$f_i(t)\geq \varphi(t)$$
for all $(i, t)\in V\times [0, T)$ such that $\varphi(t)$ exists.

Similarly, suppose $f: V\times [0, T)\rightarrow \mathbb{R}$ be a $C^1$ function such that
$$\frac{\partial f_i}{\partial t}\leq \Delta f_i+ \Phi_i(f_i), \ \ \forall (i, t)\in V\times [0, T). $$
Suppose there exists $C_2\in \mathbb{R}$ such that $f_i(0)\leq C_2$ for all $i\in V$. Let $\psi$ be the
solution to the associated ODE
\begin{equation*}
\begin{aligned}
\left\{
  \begin{array}{ll}
    \frac{d \psi}{dt}=\Phi_i(\psi), \\
    \psi(0)=C_2,
  \end{array}
\right.
\end{aligned}
\end{equation*}
then
$$f_i(t)\leq \psi(t) $$
for all $(i, t)\in V\times [0, T)$ such that $\psi(t)$ exists.
\end{theorem}

Applying the discrete maximal principle, we have the following result.
\begin{theorem}\label{converge of alpha Yamabe for negative initial curvature}
Suppose $(S, V)$ is a marked surface with a PL metric $d_0$.
$\alpha\in \mathbb{R}$ is a constant such that $\alpha \mathbf{F}_{\alpha,i}(u(0))<0$ for all $i\in V$,
then the normalized $\alpha$-Yamabe flow with surgery (\ref{alpha Yamabe flow with surgery}) exists for all time and converges
exponentially fast to a constant $\alpha$-curvature PL metric.
\end{theorem}
\proof
Note that the combinatorial $\alpha$-curvature $\mathbf{F}_{\alpha,i}$ evolves according to (\ref{evolution of curv under alpha Yamabe flow})
along the normalized $\alpha$-Yamabe flow with surgery (\ref{alpha Yamabe flow with surgery}).
The maximum principle, i.e. Theorem \ref{Maximum priciple},
is valid for this equation.
By the maximum principle, if $\alpha>0$ and $\mathbf{F}_{\alpha, i}(u(0))<0$ for all $i\in V$, we have
$$\left(\mathbf{F}_{\alpha, \min}(0)-\mathbf{F}_{\alpha, av}\right)e^{\alpha \mathbf{F}_{\alpha, av}t}\leq \mathbf{F}_{\alpha,i}-\mathbf{F}_{\alpha, av}\leq
\mathbf{F}_{\alpha, av}(1-\frac{\mathbf{F}_{\alpha, av}}{\mathbf{F}_{\alpha, \max}(0)})e^{\alpha \mathbf{F}_{\alpha, av}t}.$$
If $\alpha<0$ and $\mathbf{F}_{\alpha, i}(u(0))>0$ for all $i\in V$, we have
$$\frac{\mathbf{F}_{\alpha,av}}{\mathbf{F}_{\alpha,\min}(0)}\left(\mathbf{F}_{\alpha, \min}(0)-\mathbf{F}_{\alpha, av}\right)e^{\alpha \mathbf{F}_{\alpha, av}t}\leq \mathbf{F}_{\alpha,i}-\mathbf{F}_{\alpha, av}\leq
(\mathbf{F}_{\alpha, \max}(0)-\mathbf{F}_{\alpha, av})e^{\alpha \mathbf{F}_{\alpha, av}t}.$$
In summary, if $\alpha \mathbf{F}_{\alpha, i}(u(0))<0$ for all $i\in V$, there exists constants $C_1$ and $C_2$ such that
$$C_1e^{\alpha t\mathbf{F}_{\alpha, av}}\leq \mathbf{F}_{\alpha,i}(u(t))-\mathbf{F}_{\alpha, av}\leq C_2e^{\alpha t\mathbf{F}_{\alpha, av}},$$
which implies the long-time existence and exponential convergence of the normalized $\alpha$-Yamabe flow with surgery (\ref{alpha Yamabe flow with surgery}).
\qed

\textbf{Proof of Theorem \ref{main theorem existence of constant alpha curvature metric}:}
In the case $\alpha\chi(S)=0$, we have $\alpha=0$ or $\chi(S)=0$.
For $\alpha=0$, $\alpha$-curvature $\mathbf{F}_\alpha$ is the classical discrete curvature $\mathbf{F}$. The existence of constant curvature PL metric
is ensured by Theorem \ref{Euclidean discrete uniformization}.
If $\chi(S)=0$, the constant $\alpha$-curvature metric is a zero $\alpha$-curvature metric for all $\alpha\in \mathbb{R}$.
Especially, it is a PL metric with zero $\mathbf{F}$ curvature, the existence of which is ensured by Theorem \ref{Euclidean discrete uniformization}.

In the case of $\alpha\chi(S)<0$, by Theorem \ref{Euclidean discrete uniformization}, there is a PL metric $d'\in \mathcal{D}(d_0)$
with constant $\mathbf{F}$ curvature $\frac{2\pi\chi(S)}{N}$,
which implies that the combinatorial $\alpha$-curvature of $d'$ satisfies $\alpha \mathbf{F}_{\alpha, i}<0$ for all $i\in V$.
Applying Theorem \ref{converge of alpha Yamabe for negative initial curvature} with initial metric $d'$
gives the conclusion. \qed

\begin{remark}
Theorem \ref{main theorem existence of constant alpha curvature metric} and Theorem \ref{equivalene for alpha Yamabe flow}
together implies the long time existence and convergence of the combinatorial $\alpha$-Yamabe flow.
i.e. the $\alpha$-Yamabe flow part of Theorem \ref{main theorem convergence of alpha Yamabe and Calabi flow with surgery}.
This confirms a generalized Luo conjecture \cite{L1} on convergence of combinatorial Yamabe flow with surgery.
\end{remark}

\begin{remark}
There is another way to extend the combinatorial Yamabe flow initiated by Ge-Jiang \cite{GJ0}.
Ge-Jiang's extension comes from \cite{BPS, L2} and is designed for surfaces with fixed triangulations.
For a fixed triangulated surface, Ge-Jiang's extension of $\alpha$-Yamabe flow ensures
the long-time existence of the extended flow,
while the extended $\alpha$-Yamabe flow may converge to a virtual constant $\alpha$-curvature PL metric.
Gu-Luo-Sun-Wu's extension we use here ensures the combinatorial $\alpha$-Yamabe flow with surgery converges to a
real constant $\alpha$-curvature metric. Furthermore, Gu-Luo-Sun-Wu's extension could be applied to
extend the combinatorial $\alpha$-Calabi flow, while Ge-Jiang's extension is not valid
for this case. The readers can refer to
Subsection \ref{subsection alpha Calabi flow with surgery}
for the $\alpha$-Calabi flow with surgery.
\end{remark}

\subsection{Combinatorial $\alpha$-Calabi flow with surgery} \label{subsection alpha Calabi flow with surgery}
We can also define the combinatorial $\alpha$-Calabi flow with surgery.

\begin{definition}\label{alpha Calabi flow with surgery}
Suppose $d_0$ is a PL metric on a marked surface $(S, V)$ and $\alpha\in \mathbb{R}$.
The combinatorial $\alpha$-Calabi flow with surgery on $(S, V)$ is defined as
\begin{equation}\label{alpha Calabi flow equation with surgery}
\begin{aligned}
\left\{
  \begin{array}{ll}
    \frac{du_i}{dt}=(\Delta_\alpha \mathbf{F}_{\alpha})_i, & \hbox{ } \\
    u_i(0)=0, & \hbox{ }
  \end{array}
\right.
\end{aligned}
\end{equation}
where 
$\Delta_\alpha$ is the discrete $\alpha$-Laplace operator of $\mathbf{A}^{-1}(p, w(u(t)))\in \mathcal{D}(d_0)$
on $(S, V)$ defined by (\ref{alpha Laplacian}).
\end{definition}

Similar to the combinatorial $\alpha$-Yamabe flow on discrete Riemann surface,
$\sum_{i=1}^N w_i^\alpha$ ($\sum_{i=1}^Nu_i$ for $\alpha=0$) is invariant along the
combinatorial $\alpha$-Calabi flow with surgery (\ref{alpha Calabi flow equation with surgery}).
It is straightway to check that if the combinatorial $\alpha$-Calabi flow with surgery (\ref{alpha Calabi flow equation with surgery})
converges, the limit metric is a constant $\alpha$-curvature PL metric.

We have the following result for combinatorial $\alpha$-Calabi flow with surgery (\ref{alpha Calabi flow equation with surgery}),
which proves the combinatorial $\alpha$-Calabi flow part of Theorem \ref{main theorem convergence of alpha Yamabe and Calabi flow with surgery}.

\begin{theorem}\label{convergence of alpha Calabi flow with surgery}
Suppose $(S, V)$ is a closed connected marked surface with a PL metric $d_0$ and
$\alpha\in \mathbb{R}$ is a constant such that $\alpha\chi(S)\leq0$.
Then the combinatorial $\alpha$-Calabi flow with surgery (\ref{alpha Calabi flow equation with surgery})
exists for all time and converges to a constant $\alpha$-curvature metric in $\mathcal{D}(d_0)$.
\end{theorem}

\proof
By Theorem \ref{main theorem existence of constant alpha curvature metric},
there exists a unique PL metric $d=\mathbf{A}^{-1}(p, e^{u^*})\in \mathcal{D}(d_0)$
with $\sum_{i=1}^Ne^{\alpha u^*_i}=N$
such that $d$ has constant combinatorial $\alpha$-curvature $\mathbf{F}_{\alpha}$.

Similar to the proof of Theorem \ref{equivalene for alpha Yamabe flow}, we can define
$$W_\alpha(u)=W(u)-\mathbf{F}_{\alpha, av}\int_{u^*}^{u}\sum_{i=1}^Nw_i^\alpha du_i.$$
Then $W_\alpha$ is a well-defined $C^2$-smooth convex function defined on $\mathbb{R}^n$ under the condition $\alpha\chi(S)\leq 0$.
Furthermore, $W_\alpha(u)=W_\alpha(u+k\mathbf{1}), k\in \mathbb{R}$.
Note that $\nabla W_\alpha(u^*)=0$, we have $\lim_{u\rightarrow \infty}W_\alpha(u)|_{P}=+\infty$, where
$P=\{u\in \mathbb{R}^n|\sum_{i=1}^N w_i^\alpha=N\}$. This implies that $W_\alpha(u)|_{P}$
is a proper function on $P$.

By direct calculations, we have
\begin{equation*}
\begin{aligned}
\frac{dW_\alpha(u(t))}{dt}
=&\sum_{i=1}^N\frac{\partial W_\alpha}{\partial u_i}\frac{du_i}{dt}
=\sum_{i=1}^N(\mathbf{F}_{i}-\mathbf{F}_{\alpha, av}w_i^\alpha)(\Delta_{\alpha}\mathbf{F}_{\alpha})_i\\
=&-(\mathbf{F}_{\alpha}-\mathbf{F}_{\alpha, av})^T\cdot L\cdot(\mathbf{F}_{\alpha}-\mathbf{F}_{\alpha, av})
\leq 0.
\end{aligned}
\end{equation*}
Then $W_\alpha(u(t))$ is bounded along the combinatorial $\alpha$-Calabi flow with surgery (\ref{alpha Calabi flow equation with surgery}).
By the properness of $W_\alpha$, $u(t)$ is bounded
along the combinatorial $\alpha$-Calabi flow with surgery (\ref{alpha Calabi flow equation with surgery}),
which implies the long-time existence of combinatorial $\alpha$-Calabi flow with surgery.

As $W_\alpha(u(t))$ is bounded along the combinatorial $\alpha$-Calabi flow with surgery
and $\frac{dW_\alpha(u(t))}{dt}\leq 0$, we have
$\lim_{t\rightarrow +\infty} W_\alpha(u(t))$ exists.

Note that
\begin{equation*}
\begin{aligned}
0=&\lim_{n\rightarrow +\infty}(W_\alpha(u(n+1))-W_\alpha(u(n)))\\
=&-\lim_{n\rightarrow +\infty} (\mathbf{F}_{\alpha}-\mathbf{F}_{\alpha, av})^T\cdot L\cdot(\mathbf{F}_{\alpha}-\mathbf{F}_{\alpha, av})|_{t=\xi_n},
\end{aligned}
\end{equation*}
there is a subsequence $\xi_{n_k}$ of $\xi_n\in (n, n+1)$ such that $\mathbf{F}_\alpha(u(\xi_{n_k}))\rightarrow \mathbf{F}_{\alpha, av}$,
which implies that $u(\xi_{n_k})\rightarrow u^*\in P$.
So we have $\lim_{t\rightarrow +\infty} W_\alpha(u(t))=W_\alpha(u^*)$.
By the strictly convexity of $W_\alpha$ on $P$, we have $\lim_{t\rightarrow +\infty} u(t)=u^*$.
\qed

\begin{remark}
In the case of $\alpha=0$, the combinatorial Calabi flow with surgery was studied in \cite{ZX},
where the long time existence and convergence of the combinatorial Calabi flow with surgery were proved.
Theorem \ref{convergence of alpha Calabi flow with surgery} generalizes the results obtained in \cite{ZX}.
\end{remark}

\textbf{Acknowledgements}\\[8pt]
The authors would like to thank the referees for many valuable suggestions, which improve the paper substantially.
The research of the author is supported by Hubei Provincial Natural Science Foundation of China under grant no. 2017CFB681,
Fundamental Research Funds for the Central Universities under grant no. 2042018kf0246 and
National Natural Science Foundation of China under grant no. 61772379 and no. 11301402.

(Xu Xu) School of Mathematics and Statistics, Wuhan University, Wuhan 430072, P.R. China

E-mail: xuxu2@whu.edu.cn\\[2pt]

\end{document}